\documentclass{IEEEtran4PSCC}
\ifCLASSINFOpdf
   \usepackage[pdftex]{graphicx}
\else
   \usepackage[dvips]{graphicx}
\fi
\usepackage[cmex10]{amsmath}
\usepackage[cmex10]{amsmath}

\usepackage{cite}
\usepackage{amsmath,amssymb,amsfonts}
\usepackage{algorithmic}
\usepackage{graphicx}
\usepackage{textcomp}
\usepackage{xcolor}
\usepackage{multicol, blindtext}
\usepackage{xspace,empheq,fancybox,amsmath,amssymb,graphicx,epstopdf,epsfig,syntonly,times,amsthm} \usepackage{psfrag,color,bm,array,cite}
\usepackage{url,cite,footnote,xspace,syntonly,algorithm,algorithmic}
\usepackage{verbatim,multirow}
\usepackage[T1]{fontenc}
\usepackage{mathtools}
\usepackage{subfig}
\usepackage{amsmath,amsfonts,amsthm,bm}
\usepackage{graphicx}
\usepackage{mwe}
\usepackage{caption}
\usepackage{stfloats}
\usepackage{amsfonts}
\usepackage{stackengine}
\usepackage{lipsum}

\def\BibTeX{{\rm B\kern-.05em{\sc i\kern-.025em b}\kern-.08em
    T\kern-.1667em\lower.7ex\hbox{E}\kern-.125emX}}
    
\newcolumntype{P}[1]{>{\centering\arraybackslash}p{#1}}

\makeatletter
\let\old@ps@headings\ps@headings
\let\old@ps@IEEEtitlepagestyle\ps@IEEEtitlepagestyle
\def\psccfooter#1{%
    \def\ps@headings{%
        \old@ps@headings%
        \def\@oddfoot{\strut\hfill#1\hfill\strut}%
        \def\@evenfoot{\strut\hfill#1\hfill\strut}%
    }%
    \def\ps@IEEEtitlepagestyle{%
        \old@ps@IEEEtitlepagestyle%
        \def\@oddfoot{\strut\hfill#1\hfill\strut}%
        \def\@evenfoot{\strut\hfill#1\hfill\strut}%
    }%
    \ps@headings%
}
\makeatother

\begin{document}
%
\title{Transmission Switching Under Wind Uncertainty Using Linear Decision Rules}

\author{
\IEEEauthorblockN{Yuqi Zhou and Hao Zhu}
\IEEEauthorblockA{Department of Electrical and Computer Engineering \\
The University of Texas at Austin\\
Austin, TX, USA\\
\{zhouyuqi, haozhu\}@utexas.edu}
\and
\IEEEauthorblockN{Grani A. Hanasusanto}
\IEEEauthorblockA{Graduate Program in Operations Research \& Industrial Engineering \\
The University of Texas at Austin\\
Austin, TX, USA\\
\{grani.hanasusanto\}@utexas.edu}
}


\maketitle

\begin{abstract}
Increasing penetration of wind and renewable generation poses significant challenges to the power system operations and reliability. This paper considers the real-time optimal transmission switching (OTS) problem for determining the generation dispatch and network topology that can account for uncertain energy resources. To efficiently solve the resultant two-stage stochastic program, we propose a tractable linear decision rule (LDR) based approximation solution that can eliminate the uncertain variables and lead to fixed number of constraints. The proposed LDR approach can guarantee feasibility, and significantly reduces the computational complexity of existing approaches that grows with the number of randomly generated samples of uncertainty. Numerical studies on IEEE test cases demonstrate the high approximation accuracy of the proposed LDR solution and its computational efficiency for real-time OTS implementations. 
%
\end{abstract}


\thanksto{\protect\rule{0pt}{0mm} 
This work has been supported by NSF CAREER Grant $\#$1802319 and NSF CAREER Grant $\#$1752125.
}

\section{Introduction}
Recent years have witnessed the continued growth of wind energy and other renewable sources, which currently provide a significant portion of the electricity generation capacity in the U.S. and worldwide. The uncertainty and lack of flexibility of these renewable generation greatly challenge the real-time operational tasks for power transmission systems. One example of such tasks is the optimal transmission switching (OTS) problem by controlling the network topology to effectively reduce the operational costs and mitigate grid congestion \cite{fisher2008optimal}. With increasing variability and intermittency in energy resources, it is imperative to design efficient real-time OTS strategies that can account for problem uncertainty. 

The OTS problem under uncertainty has been considered in \cite{qiu2015chance}, which proposes a chance-constrained formulation solved by sample average approximation (SAA). To tackle the computational inefficiency due to SAA, algorithmic modifications are introduced therein, yet still requiring  uncertainty sampling. Moreover, the scenario reduction technique and heuristic methods are applied in \cite{shi2016wind} to accelerate the SAA solution. In addition, a probabilistic topology control scheme is developed in \cite{dehghanian2016probabilistic}, with
the uncertainty modeled via the point estimation to lessen the computational burden due to large-scale scenarios. To the best of our knowledge, all existing methods follow the SAA solution framework and suffer from the scalability issue. 

The goal of this work is to propose a real-time solution for OTS under uncertainty by leveraging the \textit{linear decision rules (LDR)} that can approximate the solutions of two-stage stochastic programs \cite{Chen.Sim.Peng.Zhang.2008,kuhn2011primal}. To solve the latter, one can make a \emph{here-and-now} decision at the first stage followed by adjusting \emph{recourse} actions at the second stage upon observing the uncertainty.
To eliminate the complexity arising from uncertainty sampling in SAA, the LDR approximation is an effective technique by restricting recourse decisions to be unknown affine functions of uncertainty. By assuming all uncertainty variables falling within a polytope, such affine reformulation allows to eliminate the uncertainty variables through dualization and lead to finite number of constraints. 
Hence, the main contribution of our work is three-fold:
\begin{enumerate}
  \item A \emph{two-stage} OTS problem under uncertainty is formulated which incorporates \emph{real-time} generation set-points and automatic generation control under reserve limits.
  \item A \emph{scenario-free} tractable LDR approximation is proposed for the \textit{primal} recourse decisions that can achieve high computational efficiency for the two-stage OTS problem with \textit{feasible} solutions in real time.
  \item A \emph{dual} LDR counterpart is also implemented to obtain a lower bound for the two-stage OTS problem, which together with the upper bound from primal LDR can quantify  the approximation accuracy.
\end{enumerate}

This paper is organized as follows. Section \ref{sec:ots} presents the deterministic OTS formulation and the uncertainty modeling. Section \ref{sec:PF} first formulates the two-stage stochastic program for the OTS problem under uncertainty and then develops the LDR approximations for efficient solvers. Numerical studies on the IEEE 14-bus and 118-bus systems are presented in Section \ref{sec:cs} to demonstrate the validity and efficiency of the proposed LDR approaches.




\section{System Modeling} \label{sec:ots}
We adopt the dc power flow based optimal transmission switching (OTS) formulation which does not account for reactive power or stability analysis. It is however possible to perform post-selection ac power flow analysis as in \cite{goldis2015ac}. Consider a transmission system with $N$ buses collected in the set $\cal N :=$ $\{1,\ldots,N\}$ and $L$ lines in $\cal L :=$ $\{(i,j)\} \subset \cal N \times \cal N$. For bus $i$, let $\theta_i$ be its voltage angle and collect all $\{\theta_i\}$ in  $\bm{\theta} \in \mathbb{R}^{N}$. Similarly, denote $\bm{g},~\bm{d} \in \mathbb{R}^{N}$ as vectors of generation and load per bus, respectively. Line flows $\{f_{ij}\}$ are collected in $\bm{f} \in \mathbb{R}^{L}$, given by:
\begin{align}
\bm{f} = \mathbf{K}\bm{\theta} \label{eq:PF1}
\end{align}
where matrix $\mathbf{K} \in \mathbb R^{L\times N}$ depends on the network topology. Specifically, the row corresponding to line $(i,j)$ is $b_{ij} (\mathbf{e}_i-\mathbf{e}_j)^\mathsf T$, with $b_{ij}$ being the inverse of line reactance and $\mathbf{e}_i$ as the standard basis vector. Furthermore, power conservation leads to the net injection $\bm{p} = \bm{g} - \bm{d}$ as:
\begin{align}
\bm{p} = \mathbf{A}\bm{f} \label{eq:PF2}
\end{align}
where $\mathbf{A} \in \mathbb{Z}^{N \times L}$ is the incidence matrix for the underlying graph $(\cal N, \cal L)$. The column of $\mathbf{A}$ corresponding to line $(i,j)$ equals to $(\mathbf{e}_i-\mathbf{e}_j)$.


The optimal transmission switching (OTS) problem aims to determine the optimal grid topology while minimizing the total generation cost under given load $\bm{d}$. For simplicity, linear generation cost is considered and $\bm{c} \in \mathbb{R}^{N}$ denotes the known vector of generation cost coefficients. 
We introduce the vector of binary decision variables $\bm{z} \in \mathbb{R}^{L}$ to indicate the transmission line status (1: closed, 0: open). Under a maximum number of $L_o$ open lines, the OTS problem is formulated as a mixed-integer linear program (MILP), given by
\begin{subequations} \label{eq:TS}
\begin{align}
\min \quad & {\bm{c}^{\mathsf T}\bm{g}}\\
\textrm{s.t.} \quad &  \bm{g} \in \mathbb{R}^{N}, \bm{\theta} \in \mathbb{R}^{N}, \bm{f} \in \mathbb{R}^{L}, \bm{z} \in \mathbb{Z}^{L}\label{eq:OTS_b}\\ 
  & \bm{g}^{\min} \leq \bm{g} \leq \bm{g}^{\max}    \label{eq:OTS_c}\\
  &\boldsymbol{\theta}^{\min} \leq \boldsymbol{\theta} \leq \boldsymbol{\theta}^{\max}\label{eq:OTS_d}\\
  & \bm{D}_{f}^{\min} \bm{z} \leq \bm{f} \leq \bm{D}_{f}^{\max}\bm{z} \label{eq:OTS_e}\\
  & \mathbf{A}\bm{f} = \bm{g} - \bm{d}\label{eq:OTS_f}\\
  & \mathbf{K}\boldsymbol{\theta}  - \bm{f} + \bm{D}_{\mathrm{M}}({\mathbf 1} - \bm{z}) \geq {\bm{0}}\label{eq:OTS_g}\\
  & \mathbf{K}\boldsymbol{\theta}  - \bm{f} - \bm{D}_{\mathrm{M}}({\mathbf 1} - \bm{z}) \leq {\bm{0}}\label{eq:OTS_h}\\
  & \bm{1}^{\mathsf T}\bm{z} \geq L - L_o. \label{eq:OTS_i}
\end{align}
\end{subequations}
The diagonal matrices $\bm{D}_{f}^{\min}$ and $\bm{D}_{f}^{\max}$ collect the given lower/upper limits of line flow, respectively. In addition, each diagonal entry of the diagonal matrix $\bm{D}_{\mathrm{M}}$ is a positive constant $\mathrm{M}_{ij}$ corresponding to line $(i,j)$, which will be defined soon.

Constraints in \eqref{eq:TS} are discussed in detail here. Generation, phase angle and line flow limits are set in \eqref{eq:OTS_c}-\eqref{eq:OTS_e} according to system operating limits. Note that for any open line ($z_{ij}=0$), its flow $f_{ij}$ is set to zero in \eqref{eq:OTS_e}. Constraint \eqref{eq:OTS_f} enforces network power balance in \eqref{eq:PF2}. Moreover, the pair of constraints \eqref{eq:OTS_g}-\eqref{eq:OTS_h} are introduced for the line flow model in \eqref{eq:PF1}. For any closed line ($z_{ij}=1$), the two inequalities are equivalent to an equality constraint as in \eqref{eq:PF1}. Otherwise, under $z_{ij}=0$ and thus $f_{ij}=0$ [cf. \eqref{eq:OTS_e}], these two constraints respectively become $b_{ij}(\theta_{i} - \theta_{j}) + \mathrm{M}_{ij} \geq 0$ and $b_{ij}(\theta_{i} - \theta_{j}) - \mathrm{M}_{ij} \leq 0$, both of which hold for a sufficiently large constant $\mathrm{M}_{ij}$ and are thus redundant. This is known as the \textit{Big-M} method \cite{griva2009linear}, which is popular for handling constraints with binary variables. Typically, smaller $\mathrm{M}_{ij}$ values can reduce the computation time for \eqref{eq:TS}. Thus, for each line $(i,j)$, we set: 
\begin{align}
\mathrm{M}_{ij} \coloneqq b_{ij} \Delta \theta_{ij}^{\max}
\label{eq:big_M}
\end{align}
where $\Delta \theta_{ij}^{\max}$ is a given line limit for angle stability [cf. \eqref{eq:OTS_d}]. Lastly,  \eqref{eq:OTS_i} limits the total number of open lines. 

Furthermore, we present the uncertainty model due to renewable generation. Let $K$ denote the total number of wind farms. During the window of generation dispatch or OTS, the actual wind generation $\bm{g}_w \in \mathbb{R}^{K}$ can be thought of as a deterministic nominal value $\widehat{\boldsymbol{\xi}} \in \mathbb{R}^{K}$ augmented by an uncertainty output $\boldsymbol{\xi} \in \mathbb{R}^{K}$.  Given $\boldsymbol{\xi}$ follows a distribution $\mathbb{P}$, we can represent the random vector 
\begin{align}
\bm{g}_w = \widehat{\boldsymbol{\xi}} + \boldsymbol{\xi}, \quad  \boldsymbol{\xi} \sim {\mathbb{P}}.
\label{eq:wind_stochastic}
\end{align}
For a simplified model, a more conservative assumption than distribution $\mathbb{P}$ is typically used with a given support set $\boldsymbol{\Xi}$; i.e., $\boldsymbol{\xi} \in  \boldsymbol{\Xi}$ always holds. Under a data-driven environment, the support set of the uncertainty can be easily obtained from the historical data and therefore is a realistic assumption. Instead of giving explicit realizations of the uncertainty, the conservative formulation is able to incorporate all the scenarios within the support set.
If $\boldsymbol{\Xi}$  is a closed set, then it is always possible to find upper/lower bounds for each entry. Hence, to acquire the set $\boldsymbol{\Xi}$, one can either utilize historical data or employ uncertainty bounds \cite{chen2017new} where a box constraint set can be given: 
\begin{align}
\boldsymbol{\xi}^{\min} \leq \boldsymbol{\xi} \leq \boldsymbol{\xi}^{\max}. \label{eq:wind_limit}
\end{align}
This model can be also represented using linear inequality constraints as $\mathbf{S}\boldsymbol{\xi} \leq \mathbf{t}$, with $\mathbf{S} = [ \mathbf{I};~-\mathbf{I}]$ and $\mathbf{t} = [\boldsymbol{\xi}^{\max};~-\boldsymbol{\xi}^{\min}]$. Note that this model can be used for any convex polytope constraint sets that may not be box. The ensuing section will present a fast solver for OTS problem under uncertainty with the general \emph{full-dimensional} \cite{habets2004control} uncertainty model $\mathbf{S}\boldsymbol{\xi} \leq \mathbf{t}$.

\section{OTS Problem under Uncertainty}\label{sec:PF}
To incorporate uncertainty in \eqref{eq:wind_stochastic} into OTS, we can formulate it as a two-stage stochastic program \cite{shapiro2007tutorial} by making a \emph{here-and-now} decision and taking \emph{recourse} or \emph{wait-and-see} actions once the realizations of $\boldsymbol{\xi}$ are observed. 
For the two-stage OTS problem, recourse actions on generation, phase angle, and line flow are updated once the wind output uncertainty is realized. To model the generation adjustment, we consider the automatic generation control (AGC) scheme widely used in power system operations for quickly restoring power imbalance \cite[Ch. 9]{wood2013power}. The AGC mechanism works by assigning a fixed percentage of instantaneous power imbalance to each dispatchable generator. As the total power imbalance due to $\boldsymbol{\xi}$ equals to  $\bm{1}^{\mathsf T} \boldsymbol{\xi}=\xi_1+\ldots+\xi_K$, the generation adjustment under the AGC becomes ${\bm{g}}' (\boldsymbol{\xi})=\boldsymbol{\gamma} (\bm{1}^{\mathsf T} \boldsymbol{\xi})$ with vector $\boldsymbol{\gamma}$ denoting the AGC coefficients per generator.
As for the line flows and phase angles, they are not adjusted in a fixed manner and modeled using the recourse functions ${\bm{f}}':\mathbb R^K\rightarrow \mathbb R^L$ and ${\boldsymbol{\theta}}':\mathbb R^K\rightarrow \mathbb R^N$. 
The flexibility of generator dispatch is limited by the committed reserves, with $\bm{r}^{+}$ and $\bm{r}^{-}$ denoting the upper/lower reserve limits. In addition, vector $\bm{q}$ collects the linear cost coefficients for reserve dispatch. Thus, the two-stage OTS problem under the uncertain $\boldsymbol{\xi}$ is cast as:
\begin{subequations}\label{eq:OTS_2}
\begin{align}
\min \quad &  \bm{c}^{\mathsf T}\bm{g} + \mathbb{E}[\bm{q}^{\mathsf T}\bm{g}' (\boldsymbol{\xi})] \label{eq:OTS2_a}\\
\textrm{s.t.} \quad & \eqref{eq:OTS_b} - \eqref{eq:OTS_i},\; \boldsymbol{\gamma} \in \mathbb{R}^{N} \label{eq:OTS2_b}\\
& {\bm{g}}':\mathbb R^K\rightarrow \mathbb R^N,  {\boldsymbol{\theta}}':\mathbb R^K\rightarrow \mathbb R^N, 
{\bm{f}}':\mathbb R^K\rightarrow \mathbb R^L
\label{eq:OTS2_c}\\
&\bm{g}'=\boldsymbol{\gamma} \bm{1}^{\mathsf T} \boldsymbol{\xi},\;\; \bm{r}^{-}   \leq \bm{g}'  \leq \bm{r}^{+}\;\; {\mathbb{P}} \text{-a.s.} \label{eq:OTS2_d}\\
  & \bm{g}^{\min} \leq \bm{g} + \bm{g}'  \leq \bm{g}^{\max}\;\; {\mathbb{P}} \text{-a.s.}   \label{eq:OTS2_e}\\
  & \boldsymbol{\theta}^{\min} \leq \boldsymbol{\theta} + {{{\boldsymbol\theta}}'} \leq \boldsymbol{\theta}^{\max}\;\; {\mathbb{P}} \text{-a.s.}\label{eq:OTS2_f}\\
  & \bm{D}_{f}^{\min}\bm{z} \leq \bm{f} + {{\bm{f}}'} \leq \bm{D}_{f}^{\max}\bm{z}\;\; {\mathbb{P}} \text{-a.s.}\label{eq:OTS2_g}\\
  & \mathbf{A}({\bm{f} + {{\bm{f}}'}}) = \bm{g} + \bm{g}'  - \bm{d} + \mathbf{F}(\widehat{\boldsymbol{\xi}} + {\boldsymbol\xi})\;\; {\mathbb{P}} \text{-a.s.} \label{eq:OTS2_h}\\
  & \mathbf{K}(\boldsymbol{\theta} + {{{\boldsymbol\theta}}'}) - \bm{f} - {{\bm{f}}'} + \bm{D}_{\mathrm{M}}(\bm{1} - \bm{z}) \geq \bm{0} \;\; {\mathbb{P}} \text{-a.s.}\label{eq:OTS2_i}\\
  & \mathbf{K}(\boldsymbol{\theta} + {{{\boldsymbol\theta}}'}) - \bm{f} - {{\bm{f}}'} - \bm{D}_{\mathrm{M}}(\bm{1} - \bm{z}) \leq \bm{0} \;\; {\mathbb{P}} \text{-a.s.}\label{eq:OTS2_j}
\end{align}
\end{subequations}
Note that the first-stage decisions $\{\bm{g},\bm{\theta},\bm{f},\bm{z},\boldsymbol{\gamma}\}$ also include the AGC coefficients that can be adjusted to reduce the total cost, as opposed to a fixed $\boldsymbol{\gamma}$. In addition,  the second-stage decision variables $\{{\boldsymbol{\theta}}'(\boldsymbol{\xi}), {\bm{f}}'(\boldsymbol{\xi})\}$ capture the actual adjustments on phase angle and line flow, respectively. The transformation matrix $\mathbf{F} \in \mathbb{R}^{N \times K}$ is introduced to map the wind generation  $(\widehat{\boldsymbol{\xi}} + \boldsymbol{\xi}) \in \mathbb{R}^{K}$ to the system dimension $N$.

The problem aims to minimize the sum of total generation cost at the first-stage and the expected cost of real-time \textit{recourse} reserve adjustment. Thus, constraint \eqref{eq:OTS2_d} imposes the operating reserve limits and \eqref{eq:OTS2_e}-\eqref{eq:OTS2_j} are to ensure the corresponding constraints from OTS problem \eqref{eq:TS} still hold after the recourse actions are taken [cf. \eqref{eq:OTS_c}-\eqref{eq:OTS_h}].
A two-stage stochastic program like \eqref{eq:OTS_2} is generally intractable  \cite{hanasusanto2016comment} as it involves infinitely many constraints and decision variables parameterized by $\boldsymbol\xi \in \mathbb{P}$. 
One popular solution scheme is the sample average approximation (SAA) method by approximating $\mathbb{P}$ using a finite number of scenarios $\{\boldsymbol\xi^1,\ldots,\boldsymbol\xi^S\}$ sampled  from the probability distribution $\mathbb P$. 
This way, the expectation in the objective function \eqref{eq:OTS2_a} reduces to the average ${\frac{1}{S}}\sum_{s=1}^{S} \bm{q}^{\mathsf T}\bm{g}' (\boldsymbol{\xi}^{s})$, while all the $\mathbb P$-a.s.~constraints are relaxed to finitely many constraints imposed only at the sample scenarios. The SAA method is \textit{asymptotically consistent}; i.e.,  its optimal value asymptotically converges to that of the original two-stage stochastic program as the number of samples $S$ increases. However, a very large $S$ may be necessary to gain an acceptable accuracy. This implies solving a large-scale MILP where the number of constraints and number of decision variables both grow linearly with $S$ ($\mathcal O(S)$). The computational scalability issue can prevent the SAA method from solving OTS in \textit{real time}. To address this issue, we propose a \emph{scenario-free} fast solution that utilizes the linear decision rule (LDR) scheme.

\subsection{Linear Decision Rule (LDR) Based Approximations} \label{sec:ldr}

The LDR based approximation is a powerful scheme for addressing the computational complexity of multi-stage stochastic and robust programs.
In this approximation, the recourse actions are restricted to have an affine dependence on the uncertain parameters $\boldsymbol{\xi}$, thereby using ``linear decision'' variables. Specifically, the phase angle and line flow functions are set to follow the form ${{{\boldsymbol\theta}}'}(\boldsymbol{\xi})=\bm{Y}_{\theta}{\boldsymbol\xi} + {\bm{y}}_{{\theta}}$ and ${{\bm{f}}'}(\boldsymbol{\xi})=\bm{Y}_{f}{\boldsymbol\xi} + {{\bm{y}}_{f}}$. Note that the generation adjustment ${\bm{g}}' (\boldsymbol{\xi})=\boldsymbol{\gamma} \bm{1}^{\mathsf T} \boldsymbol{\xi}$ requires no further modification as it is already in the affine form. This restriction of recourse actions leads to a conservative approximation such that its optimal value is an \textit{upper bound} on that of \eqref{eq:OTS_2}. In addition, its solution if feasible, is still feasible to \eqref{eq:OTS_2}. 
Using a given support set $\boldsymbol{\Xi}$, the LDR reformulation of \eqref{eq:OTS_2} becomes the following finite-dimensional  problem that optimizes the unknown coefficients ${\bm{Y}}_{{\theta}}$, ${\bm{y}}_{{\theta}}$,  $\bm{Y}_{f}$, and ${\bm{y}}_{f}$:
\begin{subequations}\label{eq:LDR}
\begin{align}
\min \quad & \bm{c}^{\mathsf T}\bm{g} + \bm{q}^{\mathsf T}\boldsymbol{\gamma} \bm{1}^{\mathsf T} \boldsymbol{\mu}\label{eq:LDR_a}\\
\textrm{s.t.} \quad & \eqref{eq:OTS_b} - \eqref{eq:OTS_i},\; \boldsymbol{\gamma} \in \mathbb{R}^{N} \label{eq:LDR_b}\\
  & {\bm{Y}}_{{\theta}} \in \mathbb{R}^{N \times K}, {\bm{y}}_{{\theta}} \in \mathbb{R}^{N }, \bm{Y}_{f} \in \mathbb{R}^{L \times K}, {\bm{y}}_{f} \in \mathbb{R}^{L} \label{eq:LDR_c}\\
  &\bm{r}^{-}   \leq \boldsymbol{\gamma} \bm{1}^{\mathsf T} \boldsymbol{\xi} \leq \bm{r}^{+}\;\; \forall {\boldsymbol\xi} \in  \boldsymbol{\Xi} \label{eq:LDR_d}\\
  & \bm{g}^{\min} \leq \bm{g} + \boldsymbol{\gamma} \bm{1}^{\mathsf T} \boldsymbol{\xi} \leq \bm{g}^{\max}\;\; \forall {\boldsymbol\xi} \in  \boldsymbol{\Xi} \label{eq:LDR_e}\\
  & \boldsymbol{\theta}^{\min} \leq \boldsymbol{\theta} + \bm{Y}_{\theta}{\boldsymbol\xi} + {\bm{y}}_{{\theta}} \leq \boldsymbol{\theta}^{\max}\;\; \forall {\boldsymbol\xi} \in  \boldsymbol{\Xi} \label{eq:LDR_f}\\
  & \bm{D}_{f}^{\min}\bm{z} \leq \bm{f} + \bm{Y}_{f}{\boldsymbol\xi} + {{\bm{y}}_{f}} \leq \bm{D}_{f}^{\max}\bm{z}\;\; \forall {\boldsymbol\xi} \in  \boldsymbol{\Xi}     \label{eq:LDR_g}\\
  & \mathbf{A}({\bm{f} + \bm{Y}_{f}{\boldsymbol\xi}} + {{\bm{y}}_{f}}) = \bm{g} + \boldsymbol{\gamma}\bm{1}^{\mathsf T}{\boldsymbol\xi} - \bm{d} + \mathbf{F}(\widehat{\boldsymbol{\xi}} + {\boldsymbol\xi}) \nonumber\\
  & \qquad \qquad \qquad \qquad \qquad \qquad \qquad\forall {\boldsymbol\xi} \in  \boldsymbol{\Xi}    \label{eq:LDR_h}\\
  & \mathbf{K}(\boldsymbol{\theta} + \bm{Y}_{\theta}{\boldsymbol\xi} + {\bm{y}}_{{\theta}}) - \bm{f} - (\bm{Y}_{f}{\boldsymbol\xi} + {{\bm{y}}_{f}}) \nonumber\\  & \qquad \qquad \qquad + \bm{D}_{\mathrm{M}}(\bm{1} - \bm{z}) \geq \bm{0}\;\;\forall {\boldsymbol\xi} \in  \boldsymbol{\Xi}\label{eq:LDR_i} \\
  & \mathbf{K}(\boldsymbol{\theta} + \bm{Y}_{\theta}{\boldsymbol\xi} + {\bm{y}}_{{\theta}}) - \bm{f} - (\bm{Y}_{f}{\boldsymbol\xi} + {{\bm{y}}_{f}}) \nonumber\\ & \qquad \qquad \qquad - \bm{D}_{\mathrm{M}}(\bm{1} - \bm{z}) \leq \bm{0}\;\; \forall {\boldsymbol\xi} \in  \boldsymbol{\Xi} \label{eq:LDR_j}
\end{align}
\end{subequations}
where the mean of the uncertainty $\boldsymbol{\mu}:=\mathbb E[\boldsymbol{\xi}]$.
The $\mathbb P$-a.s. constraints in \eqref{eq:OTS_2} reduce to semi-infinite constraints in \eqref{eq:LDR} as they merely involve affine functions in $\bm\xi$. 
If the support set $\bm\Xi$ is full dimensional, any semi-infinite equality of the form $\bm H\bm\xi+\bm h=\bm 0,\;\forall\bm\xi\in\bm\Xi$, holds if and only if both $\bm H=\bm 0$ and $\bm h=\bm 0$ \cite{kuhn2011primal}.
Therefore, the linear equality \eqref{eq:LDR_h} is equivalent to the following finite equalities:
\begin{equation*}
  \mathbf{A}\bm{Y}_{f} = \boldsymbol{\gamma}\bm{1}^{\mathsf T} + \mathbf{F} \;\;\textup{and}\;\;
  \mathbf{A}{\bm{f}} + \mathbf{A}{\bm{y}_f} = \bm{g}  - \bm{d} + \mathbf{F}\widehat{\boldsymbol{\xi}}.
\end{equation*}
The inequality constraints in \eqref{eq:LDR} are also amenable to tractable finite reformulations through dualization.
For instance, consider the semi-infinite constraint in \eqref{eq:LDR_f}, specifically $\boldsymbol{\theta} + \bm{Y}_{\theta}{\boldsymbol\xi} + {\bm{y}}_{{\theta}} \leq \boldsymbol{\theta}^{\max},\; \forall {\boldsymbol\xi} \in  \boldsymbol{\Xi}$. It can be reformulated to the following conditions that hold for any $i \in \cal N$:
\begin{align} \label{eq:LDR_dual_1}
    {\theta}_{i}  + {{y}}_{{\theta},{i}} - {\theta}_{i}^{\max} & \leq \min_{\mathbf{S}\boldsymbol{\xi} \leq \mathbf{t}} -\bm{e}_{i}^{\mathsf T}\bm{Y}_{\theta}{\boldsymbol\xi}
\end{align}
where $\mathbf{e}_i$ is the standard basis vector. Dualizing the right hand side of \eqref{eq:LDR_dual_1} leads to the following maximization problem:
\begin{equation}\label{eq:LDR_dual_2}
\begin{aligned}
\min_{\mathbf{S}\boldsymbol{\xi} \leq \mathbf{t}} -\bm{e}_{i}^{\mathsf T}\bm{Y}_{\theta}{\boldsymbol\xi}
 & =\max_{\boldsymbol{\pi}_{i}  \geq \bm{0}, \boldsymbol{\pi}_{i}^{\mathsf T} \mathbf{S} = \bm{e}_{i}^{\mathsf T}\bm{Y}_{\theta} } \quad  -\boldsymbol{\pi}_{i}^{\mathsf T}  \mathbf{t}.
\end{aligned}
\end{equation}
Strong linear programming duality applies because $\boldsymbol{\Xi}$ is non-empty. Collect the dual variables $\boldsymbol{\pi}_{i}^{\mathsf T}$ in $\boldsymbol{\alpha} = [\boldsymbol{\pi}_{1}^{\mathsf T};\ldots;\boldsymbol{\pi}_{N}^{\mathsf T}]$. Thus, associating \eqref{eq:LDR_dual_1} with \eqref{eq:LDR_dual_2} implies that the constraint $\boldsymbol{\theta} + \bm{Y}_{\theta}{\boldsymbol\xi} + {\bm{y}}_{{\theta}} \leq \boldsymbol{\theta}^{\max},\; \forall {\boldsymbol\xi} \in  \boldsymbol{\Xi}$ is satisfied if and only if there exists $\boldsymbol{\alpha} \geq \bm{0}$ such that $ \boldsymbol{\alpha}\mathbf{S} = \bm{Y}_{\theta} $ and $\boldsymbol{\theta}  + {\bm{y}}_{{\theta}} - \boldsymbol{\theta}^{\max}  \leq -\boldsymbol{\alpha}\mathbf{t}$. 
Introducing $\boldsymbol{\alpha}$ as a matrix of decision variables and
applying this technique to all inequality constraints in \eqref{eq:LDR}, we can obtain the full reformulation for the LDR approximation as:
\begin{subequations}\label{eq:LDR2}
\begin{align}
\min \quad & \bm{c}^{\mathsf T}\bm{g}+ \bm{q}^{\mathsf T}\boldsymbol{\gamma} \bm{1}^{\mathsf T} \boldsymbol{\mu}  \label{eq:LDR2_a}\\
\textrm{s.t.} \quad & \eqref{eq:LDR_b}, \;\eqref{eq:LDR_c} \label{eq:LDR2_b},\; \boldsymbol{\alpha}_{i}\;\; \forall i \in \{1,2,\cdots,10\}\\
  & \bm{r}^{-}  \leq -\boldsymbol{\alpha}_{1}\mathbf{t},\;
  \boldsymbol{\alpha}_{1}\mathbf{S} = -\boldsymbol{\gamma}\bm{1}^{\mathsf T} \label{eq:LDR2_c}\\
  & -\bm{r}^{+}  \leq -\boldsymbol{\alpha}_{2}\mathbf{t},\;
   \boldsymbol{\alpha}_{2}\mathbf{S} = \boldsymbol{\gamma}\bm{1}^{\mathsf T}\\ 
  & -\bm{g} + \bm{g}^{\min} \leq -\boldsymbol{\alpha}_{3}\mathbf{t},\;
    \boldsymbol{\alpha}_{3}\mathbf{S} = -\boldsymbol{\gamma}\bm{1}^{\mathsf T}\\
  & \bm{g} - \bm{g}^{\max}  \leq -\boldsymbol{\alpha}_{4}\mathbf{t},\;
   \boldsymbol{\alpha}_{4}\mathbf{S} = \boldsymbol{\gamma}\bm{1}^{\mathsf T} \label{eq:LDR2_f}\\ 
  & -\boldsymbol{\theta} - {\bm{y}}_{{\theta}} +\boldsymbol{\theta}^{\min} \leq -\boldsymbol{\alpha}_{5}\mathbf{t},\;
   \boldsymbol{\alpha}_{5}\mathbf{S} = -\bm{Y}_{\theta}\\
  & \boldsymbol{\theta} + {\bm{y}}_{{\theta}} -\boldsymbol{\theta}^{\max} \leq -\boldsymbol{\alpha}_{6}\mathbf{t},\;
   \boldsymbol{\alpha}_{6}\mathbf{S} = \bm{Y}_{\theta}\\ 
  & \bm{D}_{f}^{\min}\bm{z} - \bm{f} - {\bm{y}_f} \leq -\boldsymbol{\alpha}_{7}\mathbf{t},\;
    \boldsymbol{\alpha}_{7}\mathbf{S} = -\bm{Y}_{f}\\ 
  & -\bm{D}_{f}^{\max}\bm{z} + \bm{f} + {\bm{y}_f} \leq -\boldsymbol{\alpha}_{8}\mathbf{t},\;
   \boldsymbol{\alpha}_{8}\mathbf{S} = \bm{Y}_{f}\\ 
  & \mathbf{A}\bm{Y}_{f} = \boldsymbol{\gamma}\bm{1}^{\mathsf T} + \mathbf{F},\;
  \mathbf{A}{\bm{f}} + \mathbf{A}{\bm{y}_f} = \bm{g}  - \bm{d} + \mathbf{F}\widehat{\boldsymbol{\xi}}\\ 
  & -\mathbf{K}(\boldsymbol{\theta}+{\bm{y}}_{{\theta}}) + \bm{f} + {\bm{y}_f} \leq \bm{D}_{\mathrm{M}}(\bm{1} - \bm{z})-\boldsymbol{\alpha}_{9}\mathbf{t}\\
  & \mathbf{K}(\boldsymbol{\theta}+{\bm{y}}_{{\theta}}) - \bm{f} - {\bm{y}_f}\leq \bm{D}_{\mathrm{M}}(\bm{1} - \bm{z})-\boldsymbol{\alpha}_{10}\mathbf{t}\\
  & \mathbf{K}\bm{Y}_{\theta} + \boldsymbol{\alpha}_{9}\mathbf{S} = \bm{Y}_{f},\;
  -\mathbf{K}\bm{Y}_{\theta} + \boldsymbol{\alpha}_{10}\mathbf{S} = -\bm{Y}_{f}.
\end{align}
\end{subequations}
Thanks to the dualization, the LDR formulation \eqref{eq:LDR2} benefits significantly from a fixed number of decision variables and constraints, as compared to the growing number $\mathcal O(S)$ in SAA. This ensures an affordable computational complexity that is very attractive for real-time operations.

To assess the suboptimality of the LDR approximation \eqref{eq:LDR2}, we develop the \textit{dual LDR} counterpart proposed in \cite{kuhn2011primal} that can obtain the \textit{lower bound} of \eqref{eq:OTS_2}. The idea is to apply the LDR on the Lagrangian multipliers for \eqref{eq:OTS2_d}-\eqref{eq:OTS2_j}, instead of the primal variables therein, by restricting the former to be affine functions of $\bm\xi$. The aforementioned technique of converting semi-infinite constraints to tractable finite-dimensional ones can be similarly applied as well. Since this approach approximates the dual problem, instead of the primal, its solution provides a lower bound for the optimal  objective of \eqref{eq:OTS_2} due to duality theory. To sum up, the primal LDR formulation in \eqref{eq:LDR2} guarantees feasibility and upper bound of the original problem \eqref{eq:OTS_2}, while the dual LDR obtains the lower bound with no feasibility guarantee. As in any approximation bounds, the gap between the two bounds is very useful for evaluating the effectiveness of the feasible solution achieved by the primal LDR \eqref{eq:LDR2}. This is especially attractive for large-scale systems where the two-stage stochastic program \eqref{eq:OTS_2} is intractable.



\section{Numerical Results} \label{sec:cs}
This section presents the numerical performance comparisons of the SAA method and the proposed LDR methods using the IEEE 14-bus and 118-bus test cases. The test cases have been slightly modified to include uncertain wind generation. Since the uncertainty bounds and the resultant support set can be obtained offline from historical data, the computation time recorded here only includes the solution time for the optimization problems.
The optimization models are implemented in the MATLAB\textsuperscript{\textregistered} and further solved by CPLEX solver on a regular laptop with Intel\textsuperscript{\textregistered} CPU @ 2.60 GHz and 12 GB of RAM.


\subsection{14-Bus System Tests}

The original IEEE 14-bus system consists of 20 lines and 5 conventional generators. We add 5 wind farms to the case, located at buses 3, 5, 6, 10, and 13, respectively. To quantify the wind output uncertainty, we use the upper/lower bound model in \eqref{eq:wind_limit}, with the range proportional to the nominal wind output $\widehat{\boldsymbol{\xi}}$ and centered at zero ($\mathbb E[\boldsymbol{\xi}]=\bm{0}$).
To test the SAA, we randomly generate samples of $\boldsymbol{\xi}$ for a sufficiently large sample size $S=500$. Although different types of probability distributions have been examined for the SAA, we use the uniform distribution within the support set as an example. Notice that one can also incorporate spatial and temporal correlation of wind \cite{xie2013short} in generating the samples, in this work for simplicity we treat the wind uncertainty at different wind farms $\xi_{i}, \forall i \in \{1,\cdots,K\}$ as independent random variables. 

Since SAA can provide high accuracy in approximating the two-stage OTS problem \eqref{eq:OTS_2}, its result is used as the benchmark solution for evaluating LDR-based methods.
\begin{table}[t!] 
\caption{Line Switching Decisions for the 14-Bus System} \label{tab:sd_14} 
\centering 
{\renewcommand{\arraystretch}{1.2}
\begin{tabular}{|P{0.8cm}|p{1.7cm}|p{1.7cm}|p{1.7cm}|}
\hline 
 $L_o$ &  SAA  & Primal LDR  &  Dual LDR \\ 
\hline 
1 & [20] & [20] & [13]\\   
2 & [19;20]  & [19;20] & [13]\\  
3 & [9;18;20]  & [9;18;20] & [9,11,20]\\  
4 & [9;13;18;20]  & [9;13;18;20] & [9,11,20]\\  
\hline  
\end{tabular}}
\end{table}
Table \ref{tab:sd_14} lists the switching decisions obtained by SAA and LDR methods under various numbers of open lines. The primal LDR gives exactly the same optimal $\bm{z}^{*}$ as SAA, while the dual LDR leads to suboptimal solutions. This observation confirms our analysis on the LDR solutions.
To quantify the accuracy of primal LDR approximation, the \emph{out-of-sample} costs are listed in Table \ref{tab:oc_14}. The out of sample test refers to re-solving the second-stage optimization problem for a new set of samples with the optimal first-stage solutions fixed. This approach can realistically depict the implementation performance for each method. Since SAA has been solved only for the original samples, its first-stage decisions may become infeasible for new scenarios. Such issue, if any, can be addressed by having larger sample size to tighten the feasible region the first time we run SAA. As the primal LDR works for any $\boldsymbol{\xi}$ in the support set, there is no infeasibility issue.
The dual LDR is neglected here as it provides the lower bound solution. Although SAA and primal LDR share the same $\bm{z}^{*}$ as shown in Table \ref{tab:sd_14}, their generation dispatch solutions  $\bm{g}^*$ slightly differ, and so do their out-of-sample costs. 
Interestingly, the out-of-sample cost for the primal LDR closely approaches that of SAA, only at 0.9$\%$ increase on average, so the generation dispatch obtained by the LDR is also a competitive approximation.

\begin{table}[t!] 
\caption{Out-of-sample Costs for SAA and LDR} 
\centering 
{\renewcommand{\arraystretch}{1.2}
\begin{tabular}
{ |c | c | c | c | c | c | c |} 
\hline 
 $L_o$ & SAA $(\$/h)$  & Primal LDR $(\$/h)$ & Deviation\\ 
\hline 
1 & 378.7  & 383.9  & 1.37$\%$\\   
2 & 377.9  & 383.0  & 1.35$\%$\\  
3 & 364.2  & 366.4  & 0.60$\%$\\  
4 & 364.1  & 365.5  & 0.38$\%$\\  
\hline  
\end{tabular} }
\label{tab:oc_14} 
\end{table}

\begin{figure}[t!]
\centering
\vspace{-2pt}
\includegraphics[trim=0cm 0cm 0cm 1cm,clip=true,totalheight=0.18\textheight]{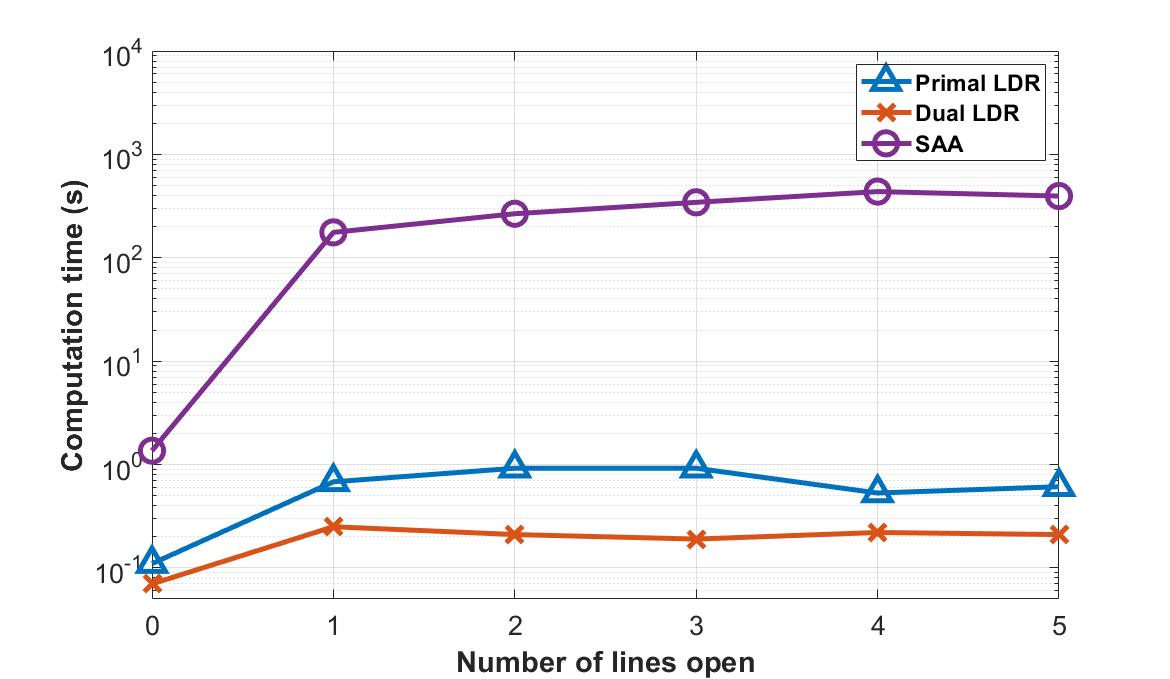}
\caption{Computation time comparison for 14-bus system.}
\label{fig:ct_14}
\end{figure}
Lastly, the computation time of all methods are plotted in Fig. \ref{fig:ct_14}. The LDR methods enjoy high computation efficiency, both solved within 1 second for various $L_o$ values. As the SAA's complexity increases with the number of samples, it takes much longer time (>$100$ seconds) to achieve reasonable accuracy as compared to LDR. This confirms the improvement of the proposed LDR solutions for real-time operations.

\subsection{118-Bus System Tests}
The larger-sized 118-bus system is also tested, consisting of 186 lines and 19 conventional generators. Five or more wind farms have been added in a similar fashion. Due to the scalability issue, the SAA method requires much higher computational time for this large system.
The solver cannot obtain a converged solution for the SAA method within an hour, even for $L_o=1$ and a reduced sample size ($S=50$).  Therefore, only the proposed LDR methods are presented here, with the approximation accuracy evaluated using the gap between the upper and lower bounds that are obtained respectively from primal and dual LDR.

\begin{figure}[t!]
\centering
\vspace{-2pt}
\includegraphics[trim=0cm 0cm 0cm 1.2cm,clip=true,totalheight=0.20\textheight]{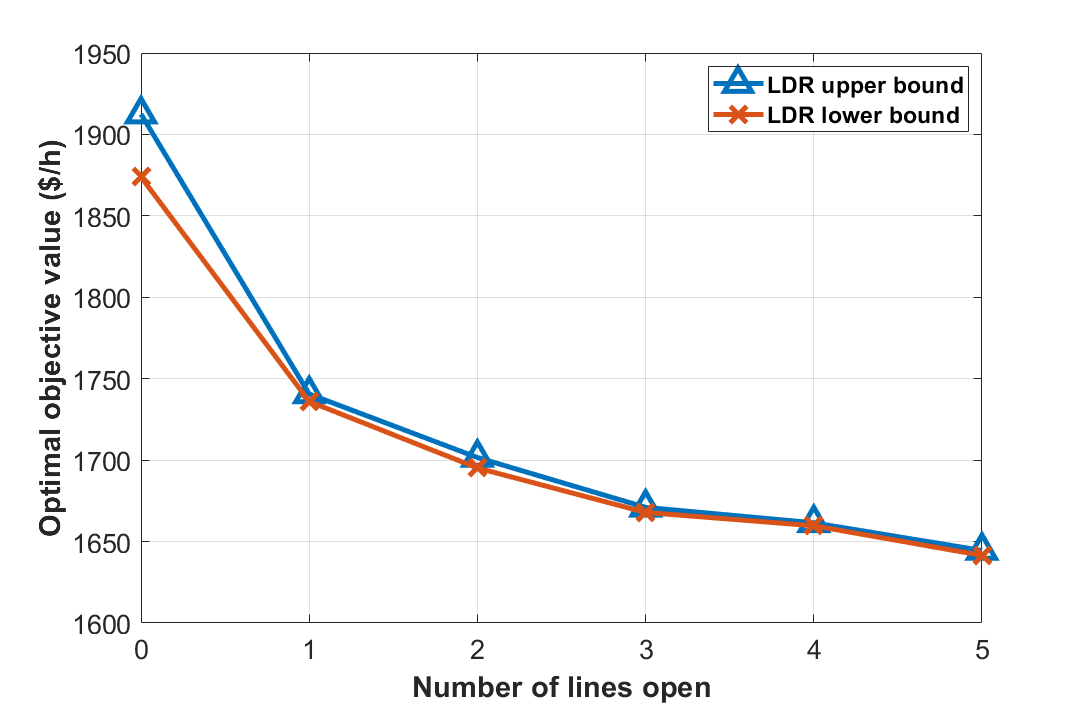}
\caption{LDR upper/lower bounds for the 118-bus system.}
\label{fig:bound_118}
\end{figure}

We first fix the number of wind farms to be $K=5$ and compare the LDR upper and lower bounds by varying $L_o$, as shown in Fig. \ref{fig:bound_118}. Overall, the gap between the two bounds is very small ($\leq 0.6\%$), validating the near-optimal performance of the primal LDR in approximating the original problem in this case. 
Furthermore, we increase the number of wind farms ($K=7,8,9,10$) and list the computation time of primal LDR in Table \ref{tab:ct_118}. The results again demonstrate the computational efficiency of the proposed LDR approach and support its real-time implementation for the OTS problem under uncertainty. 

\section{Conclusion}
\label{sec:con}
In this paper we present a transmission switching strategy under uncertainty that accommodates the real time dispatch of generators, which is formulated as a two-stage stochastic program. To tackle the tractability issue and allow real time application, we further implement the linear decision rule approximation by restricting recourse actions to be affinely dependent on the uncertainty. We show that the uncertainty can be conveniently incorporated into a scenario-free reformulation through dualization. 
Numerical studies demonstrate the accuracy of the linear decision rule approach and its significant improvement in computational efficiency.

\begin{table}[t!] 
\caption{Primal LDR run time (in sec) on 118-bus system} 
\centering 
{\renewcommand{\arraystretch}{1.2}
\begin{tabular}{|P{0.8cm}|P{1.3cm}|P{1.3cm}|P{1.3cm}|P{1.3cm}|}
\hline 
 $L_o$ & $K=7$  & $K=8$   & $K=9$ & $K=10$ \\ 
\hline 
1 & 0.31  & 0.33  & 0.34  & 0.45\\   
2 & 0.51  & 0.48  & 0.52  & 0.55\\  
3 & 1.39  & 1.63  & 1.38  & 1.45\\  
4 & 3.07  & 3.21  & 3.25  & 3.28\\  
5 & 6.70  & 6.55  & 7.13  & 7.15\\ 
\hline  
\end{tabular} }
\label{tab:ct_118} 
\end{table}

\bibliography{bibliography.bib}

\bibliographystyle{IEEEtran}

\itemsep2pt

\end{document}